\pdfoutput=1

\documentclass{article}
\usepackage[T1]{fontenc}
\usepackage[UKenglish]{babel}

\usepackage{Factory, Math, Theorema, Classic, Styl}

\usepackage[all, 2cell]{xy} \UseAllTwocells 
\newcommand{\upar}{ \ar@<+.4ex> }			\newcommand{\downar}{ \ar@<-.4ex> }

	\newcommand{\Mod}{\mathfrak{Mod}} \newcommand{\XMod}{\mathfrak{XMod}}
\newcommand{\Func}{\mathfrak{Fun}}		\newcommand{\Pol}{\mathfrak{Pol}} 	
\newcommand{\Num}{\mathfrak{Num}}				
\newcommand{\SPol}{\mathfrak{SPol}} 				\newcommand{\HPol}{\mathfrak{Hom}}
\newcommand{\Laby}{\mathfrak{Laby}}
\DeclareMathOperator{\Lin}{\mathrm{Lin}}
\DeclareMathOperator{\Ann}{\mathrm{Ann}}

\newcommand{\B}{\mathbf{B}}
\newcommand{\BB}{\B[\B^{n\times n}]_n}				
\newcommand{\GB}{\Gamma^n(\B^{n\times n})} 
\newcommand{\GH}[1][-]{\Gamma^{n}\Hom(\B^n,#1)}
\newcommand{\de}{\diamond} 			\DeclareMathOperator*{\De}{\lozenge} 
\DeclareMathOperator{\ce}{\mathrm{cr}}	
	
\newcommand{\maze}[1]{\begin{bmatrix}\xymatrixcolsep{.8pc}\xymatrixrowsep{.5pc}\xymatrix{#1}\end{bmatrix}}

\usepackage{fancyhdr}						
\begin{document}

\lefthyphenmin=2 \righthyphenmin=2

\titul{ON EXTENSIONS OF \\ POLYNOMIAL FUNCTORS}
\auctor{Qimh Richey Xantcha\footnote{\textsc{Qimh Richey Xantcha}, Uppsala University: \texttt{qimh@math.uu.se}}}
\datum{\today}

\maketitle

\bigskip 

\epigraph{\begin{vverse}[What will the Line stretch out to'th'cracke of Doome?]
What will the Line%
\footnote{``It is reported that Voltaire often laughs at the tragedy of \emph{Macbeth}, for 
having a legion of ghosts in it. One should imagine he either had not learned English, 
or had forgot his Latin; for the spirits of Banquo's line are no more ghosts, than the 
representations of the Julian race in the \AE neid; and there is no ghost but Banquo's
throughout the play.'' --- Mrs.~Elizabeth Montague: \emph{Essay on the Genius and Writings of Shakspeare} (1769)}
stretch out to'th'cracke of Doome? \\
Another yet? a seauenth? Ile see no more: \\
And yet the eight appeares, who beares a glasse, \\
Which shewes me many more: and some I see, \\
That two-fold Balles, and trebble Scepters carry. \par 
\vattr{Shakspeare, \emph{The Tragedy of Macbeth}}
\end{vverse}}

\bigskip 

\begin{argument} \noindent
A formula for calculating Extensions of (mainly integral) Polynomial Functors is established, based upon projective 
resolutions. Sample computations are performed, which, in particular, exhibit a surprising 
non-trivial extension of Divided Cubes by Symmetric Cubes. 
An explicit description of the latter is given. 
\end{argument}

\bigskip 

\noindent
The classification of Extensions of Polynomial Functors has, for several decades, 
been known to be an arduous quest. 
Significant progress has been made over fields; 
less so for more general rings, epitomised by the integers $\Z$. 

Our modest purpose shall be to gather some knowledge of the groups $\Ext^1(F,G)$, 
where $F$ and $G$ are 
chosen from among the classical algebraic functors $S^n,\Lambda^n,\Gamma^n$ (possibly of different degrees). 
Several variations upon this theme suggest themselves, depending upon the choice of base ring
as well as in which category the extensions be permitted to dwell. 
Let us survey the known results from the literature. 

We denote by $\Func$ the category of functors from 
finitely generated, free modules to modules, by $\SPol$ the strict polynomial functors, and by 
$\Pol_n$ the polynomial functors of degree (at most) $n$. 
 
\emph{1. Extensions of Strict Polynomial Functors over $\Z$.} 
\btab{|c|ccc|}
\hline
$\Ext^1_{\SPol}(-,-)$	& $S^n$ 		& $\Lambda^n$ 								& $\Gamma^n$ \\
\hline 
$S^n$					& $0$		& $\bca 0 & (n=0,1) \\ \Z_2 & (n\geq 2) \eca$	& $0$   \\
$\Lambda^n$				& $0$		& $0$										& $\bca 0 & n=0,1 \\ \Z_2 & (n\geq 2) \eca	$ \\
$\Gamma^n$ 				& $0$		 											& $0$								& $0$ \\ 
\hline
\etabc{Extensions of Strict Polynomial Functors over $\Z$. \label{T: Bar}}
The full extension groups, including, when applicable, the Yoneda ring structure, 
were described by Touzé (\cite{Bar}) in this case. 
We have reproduced in Table \ref{T: Bar} the values of $\Ext^1$. 
Not all of these calculations 
were original; work in this direction seems to have been pioneered by Akin (\cite{Akin}), who disposed of the case 
$\Ext^1_{\SPol}(\Lambda^n,\Gamma^n)=\Z_2$.

The two non-split extensions in Table \ref{T: Bar} are 
$$ 
\xymatrix{0\ar[r] & \Gamma^n \ar[r]^\alpha & Y^n \ar[r]^\beta & \Lambda^n \ar[r] & 0   } 
$$
and
$$ 
\xymatrix{0\ar[r] & \Lambda^n \ar[r]^\gamma & Y^n \ar[r]^\delta & S^n \ar[r] & 0,   } 
$$
where $Y^n$ represents the quotient of $T^n$ by the action of the alternating group. The maps $\alpha$, $\beta$, and $\delta$ are 
respectively induced by the canonical maps 
$$\Gamma^n\to T^n,\qquad T^n\to \Lambda^n,\qquad T^n\to S^n. $$
Writing $x\times y$ for the image of $x\otimes y$ in $Y(M)$, the map $\gamma$ is given by the formula 
$$ 
x_1 \wedge x_2 \wedge x_3 \wedge \cdots \wedge x_n 
\mapsto [()-(12)] (x_1 \times x_2 \times x_3 \times \cdots \times x_n) .
$$

\emph{2. Extensions of Polynomial Functors over $\Q$.} The functor category is 
then semi-simple, and all exact sequences split.

\emph{3. Extensions of Polynomial Functors over Finite Fields.}
In their study \cite{FLS} of Mac Lane cohomology, of which they have demonstrated the great value, 
Franjou, Lannes \& Schwartz successfully 
computed $\Ext_{\Func}(I,S^n)$
over finite prime fields and $\Ext_{\Func}(I,I)$ over arbitrary finite fields. 
Franjou's paper \cite{Puissances} threw some further light upon
 the extension algebras of symmetric and exterior powers over finite prime fields. 
Various $\Ext$-algebras of classical algebraic functors were then impressively 
mapped out by Franjou, Friedlander, Scorichenko \& Suslin in \cite{FFSS} (Theorem 6.3).

\emph{4. Extensions of Polynomial Functors over $\Z$.}  
\btab{|c|ccc|}
\hline
$\Ext^1_{\Func}(-,-)$			& $S^2$ 			& $\Lambda^2$ 	& $\Gamma^2$ \\
\hline 
$S^2$						& $0$			& $\Z_2$			& $0$ \\
$\Lambda^2$					& $0$			& $0$			& $\Z_2$ \\
$\Gamma^2$ 				 	& $0$			& $0$			& $0$ \\ 
\hline
\etabc{Extensions of Quadratic Functors over $\Z$. \label{T: Splitting}}
This, the least explored and, as it seems, most perplexing case, shall be our main line of enquiry.
We have located a single source \cite{Integers} knowledgeable on the subject, 
wherein Franjou and Pirashvili perform a quick and elementary calculation of the 
Yoneda ring $\Ext_{\Func}(I,I)$, as 
well as its graded modules  
$\Ext_{\Func}(I,S^n)$ and $\Ext_{\Func}(I,\Lambda^n)$.

Not withstanding this scarcity of written information,  
extensions of quadratic functors, collected in Table \ref{T: Splitting}, appear well known to the cognoscenti. 
They may be calculated using the methods advocated, for example, by Mikhailov in \cite{Splitting}.

\emph{5. Comparison of Extensions.} Consider two strict polynomial functors $F$ and $G$ of degree $n$ over a finite field $K$. 
The existence of a canonical embedding
$$ 
\Ext_{\SPol}(F,G) \subseteq \Ext_{\Func}(F,G),
$$
preserving homological degree,
was exhibited by Franjou, Friedlander, Scorichenko, and Suslin in \cite{FFSS} (Corollary 3.7).
They proceed to derive their celebrated 
\emph{Strong Comparison Theorem} (Theorem 3.10): 
$$
\indlim_{j\to\infty} \Ext_{\SPol}(F^{(j)},G^{(j)}) \cong  \Ext_{\Func}(F,G), 
$$
assuming merely that $n\leq \abs{K}$. The symbol $H^{(1)}$ represents the Frobenius 
twist of the functor $H$, and the formula thus explains the great interest taken in the extension groups of such. 

We now present our main results, 
for simplicity stated for polynomial functors over the ring $\Z$, 
although they will 
apply in the somewhat more general context of numerical functors over a binomial ring.

\begin{inttheorem}[\ref{S: ExtFP}]
When $F$ and $G$ are strict polynomial integral functors of degree $n$, 
with $G$ being torsion-free, there is an embedding
$$ \Ext^1_{\SPol}(F,G) \subseteq \Ext^1_{\Pol_n}(F,G). $$
\end{inttheorem}

It is well known (see Theorem \ref{S: ExtF}) that the category $\Pol_n$ is extension-closed, so that 
$$ \Ext^1_{\Pol_n}(F,G) = \Ext^1_{\Func}(F,G)$$
for polynomial functors $F$ and $G$. Example \ref{E: Num} below shews that the category 
$\Num_n$ of \emph{numerical} functors will not, in general, be extension-closed. 

\begin{inttheorem}[\ref{S: D}] When $F$ and $G$ are classical algebraic functors over the integers, 
each a symmetric, exterior, or divided power, 
there is an explicitly given integral matrix $D$ such that 
$ \Ext^1_{\Func}(F,G)$ is isomorphic with the torsion subgroup of $\Coker D$.
\end{inttheorem}

We next embark upon the project of applying the theorem to some concrete situations.  
The few $\Ext^1$-groups already known merit some consideration, and the results produced by our method 
will be found to comply with those from the scholarly literature.  

It will then naturally be asked what hitherto unknown $\Ext^1$-groups may be uncovered by our method.
This is the question we turn to next, with some rather striking results. 
Quite unexpectedly, we shall find in Examples \ref{E: 3} and \ref{E: 4} that
$$ \Ext^1(\Gamma^3,S^3) = \Z_2 \qquad \text{and} \qquad \Ext^1(\Gamma^4,S^4)=\Z_2 \oplus \Z_3. $$
The mysterious non-trivial extension $X$ of $\Gamma^3$ by $S^3$ deserves deliberation. 
Its definition is deceptively simple: 
$$
X(M) = \gen{p,\ x^{[3]},\ x^{[1]}y^{[1]}z^{[1]},\ x^{[2]}y^{[1]}+\frac{x^2y+xy^2}{2} | p\in S^3(M),\ x,y,z\in M},
$$
 a subfunctor of $\left(\gen{\frac12}\otimes S^3\right)\oplus \Gamma^3$. 
Full details on this functor may be found in Example \ref{E: X} below, 
 along with an elementary argument as to why its associated exact sequence does not split. 

The last section of the text will be consecrated to a systematic (but incomplete) enquiry into $\Ext^1_{\Func}(F,G)$, for 
$F\in\{S^m,\Lambda^m,\Gamma^m\}$ and $G\in \{S^n,\Lambda^n,\Gamma^n\}$, reaching the conclusion: 

\btab{|c|ccc|}
\hline
$\Ext^1_{\Func}(-,-)$			& $S^n$ 			& $\Lambda^n$ 	& $\Gamma^n$ \\
\hline 
$S^m$						& $0$ ($m\geq n$)			
											& $\bca
											\Z_2 & \text{($m= n\geq 2$)} \\
											0 & \text{($\neg m=n\geq 2$)}
												\eca$
															& $0$ \\
$\Lambda^m$					& $0$ ($m\geq n$)	& $0$			& $\bca
																\Z_2 & \text{($m= n\geq 2$)} \\
																0 & \text{($\neg m=n\geq 2$)}
																\eca $ \\
$\Gamma^m$ 				 	& ---				& $0$ ($m\leq n$)	& $0$ ($m=n$) \\ 
\hline
\etabc{Extensions of Polynomial Functors over $\Z$. \label{T: Ext}}

\begin{inttheorem}[\ref{S: Ext}]
Inside the functor category $\Func$ over the base ring $\Z$, extension groups are as indicated in Table \ref{T: Ext}.
\end{inttheorem}

The research was carried out at
the Institut de Recherche Mathématique Avancée at Strasbourg. We gratefully acknowledge its hospitality, and 
the cordial invitation received from Dr Christine Vespa. We are much obliged toward  
the Swedish--French Foundation for its generous financial support. We owe much to the helpful 
criticism of Dr Antoine Touzé.

\setcounter{section}{-1}
\section{Polynomial Functors}			\label{A: Preliminaries}

Throughout this paper, we shall work with a fixed ring $\B$ of scalars. 
All modules, homomorphisms, and tensor products will be taken over this $\B$, unless otherwise stated. 
We let $\Mod={}_\B \Mod$ betoken the category of (unital) modules over this ring, and we denote
by $\XMod$ the subcategory of finitely generated and free modules.

When $A$ and $B$ are linear categories (enriched over $\Mod$), 
the symbols $\Fun(A,B)$ and $\Lin(A,B)$ denote the 
categories of functors and linear functors, respectively, from $A$ to $B$.
A \textbf{module functor} is a functor $\XMod\to\Mod$. These functors form an abelian category denoted 
by $\Func$.

Polynomial functors were introduced by Eilenberg \& Mac Lane, \cite{EM}, Sections 8 and 9.
Let $\phi\colon M\to N$ be a map of modules. The $n$'th 
\textbf{deviation} of $\phi$ is the map 
$$ 
\phi(x_1\de\cdots\de x_{n+1}) = \sum_{I\subseteq[n+1]} (-1)^{n+1-|I|} \phi\left(\sum_{i\in I} x_i\right)
$$
of $n+1$ variables. 

\begin{definition}
The module functor $F\colon \XMod \to \Mod$ is \textbf{polynomial} of degree (at most) $n$ 
if 
$$
F(\alpha_1\de \cdots \de \alpha_{n+1}) = 0 
$$  
for any homomorphisms $\alpha_1,\dots,\alpha_{n+1}\colon M\to N$.
\end{definition}


\begin{definition}
Let $F$ be a module functor and let $k$ be a natural number. 
The \textbf{cross-effect} of rank $k$ is the multi-functor 
$$ \ce_k F(M_1,\dots,M_k) = \Im F\left(\pi_1\de\cdots\de\pi_k \right), $$
where $\pi_i\colon \bigoplus M_i \to \bigoplus M_i$ denote the canonical projections. 
\end{definition}

A functor is polynomial of degree $n$ if and only if the cross-effects of rank exceeding $n$ vanish.

For most of this paper, we shall let $\B$ be a \emph{binomial ring}%
\footnote{This is a \emph{numerical ring} in the 
terminology \cite{TE} of Ekedahl. 
As to the equivalence of the two notions, a proof is offered in \cite{BR}.} 
in the sense of Hall (\cite{Hall}); 
that is, commutative, unital, and in the possession of binomial co-efficients. 
 In particular, the reader may assume $\B=\Z$. 

For such a ring, one may speak of \emph{numerical} maps and functors, 
as introduced in \cite{PM} and \cite{PF}. 

\begin{refdefinition}[\cite{PF}, Definition 10] Assume $\B$ binomial.
The module functor $F\colon \XMod \to \Mod$ is  \textbf{numerical} of degree (at most) $n$ if 
it satisfies the following two equations: 
\begin{gather*}
F(\alpha_1\de\cdots\de \alpha_{n+1}) =0, \qquad \alpha_1,\dots,\alpha_{n+1}\colon M\to N; \\
F(r\alpha) = \sum_{k=0}^n \binom{r}{k}\phi\left(\De_k \alpha\right), \qquad r\in\B,\ \alpha\colon M\to N.
\end{gather*}
\end{refdefinition}

Over $\Z$, numerical and polynomial co-incide, but in the general case there is  
 a trade-off betwixt the two. 
On the one hand, polynomial functors behave 
well under extensions (cf.~Theorem \ref{S: ExtF} below), essentially due to 
the cross-effects being functorial in nature,
whereas numerical ones do not (cf.~Example \ref{E: Num}). 
But the structure of polynomial functors over an arbitrary base ring is not easily 
deduced. By contrast, 
numerical functors are encoded by a quotient labyrinth category whose basis is the pure mazes 
(see Theorem \ref{S: Num} below), which means their structure theory closely resembles that 
of polynomial functors over $\Z$.

Also, numerical functors allow for an easy comparison with the strict polynomial functors 
of Friedlander \& Suslin, which we now introduce. 
A striking illustration would be the Polynomial Functor Theorem of \cite{PF}.

\begin{refdefinition}[\cite{FS}, Definition 2.1] Assume $\B$ commutative and unital. 
 The module functor $F\colon \XMod \to \Mod$ is 
 \textbf{strict polynomial} of degree $n$ if the arrow maps
$$ F\colon -\otimes\Hom(M,N) \to -\otimes\Hom(F(M),F(N)) $$
are natural transformations from commutative, unital algebras to sets. 
\end{refdefinition}

Over a $\Q$-algebra, numerical and strict polynomial functors co-incide. 

Strict polynomial functors decompose as the direct sum of their homogeneous components. 
We shall denote by $\HPol_n$ the abelian category of homogeneous functors of degree $n$.

Let us now, for the moment, assume $\B$ unital only, and briefly recall 
the definition of the Labyrinth Category over $\B$, as introduced in \cite{Fun}, Definition~22; 
though its first appearance can be traced back to \cite{X}. 
Its objects are finite sets, and its arrows are 
linear combinations of so-called mazes. A \emph{maze} from the set
$X$ to the set $Y$ is a multi-set of \emph{passages}, by which is meant a formal arrow from an element
of $X$ to an element of $Y$, decorated with a scalar. For example, the following diagram 
constitutes a maze from $\{a,b,c,d\}$ to $\{w,x,y,z\}$:
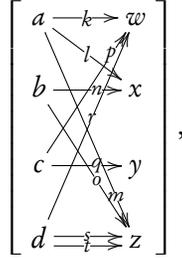
$$ 
\begin{bmatrix} 
\xymatrixrowsep{1.5em} 
\xymatrix{
a \ar[r]|{k} \ar[dr]|{l} \ar[dddr]|(0.8){m} 	& w \\
b \ar[r]|(0.6){n} \ar[ddr]|(0.6){o}	    	& x \\
c \ar[uur]|(0.75){p} \ar[r]|(0.6){q} 	    	& y \\
d \ar[uuur]|(0.55){r} \upar[r]|{s} \downar[r]|{t}& z  
} 
\end{bmatrix}, 
$$
where $k,l,m,n,o,p,q,r,s,t$ are scalars. Composition of mazes is formed 
according to a complicated rule (\cite{Fun}, Definition 21), and the 
resulting category is spoken of as 
the \emph{Labyrinth Category} $\Laby$. 
 
Being a linear category, it possesses most qualities one would normally ascribe to an associative 
algebra, except a global multiplicative identity is lacking 
since the category has infinitely many non-isomorphic objects. In particular, one may 
investigate its module theory, and the main theorem 
of \cite{Fun} is that functors of $\B$-modules are simply modules over this category.

In order to state the precise result, we need some notation. 
Let $$ \sigma_{yx}\colon\B^{X}\to\B^{Y}, \quad x\in X,\ y\in Y$$ denote the canonical \emph{transportations}, 
mapping a $1$ in position $x$ to a $1$ in position $y$, and everything else to $0$.
Moreover, let $\ce_X F(\B)$ denote the cross-effect of rank $X$ of the functor $F$, evaluated on $|X|$ copies 
of the module $\B$. 

\begin{reftheorem}[\cite{Fun}, Theorem 14] 
The functor \label{T: LoF}
$$ \Phi\colon \Func\to \Lin(\Laby,\Mod), $$
where $\Phi(F)\colon \Laby\to\Mod$ takes
\begin{gather*}
X \mapsto \ce_X F(\B) \\
[P\colon X\to Y] \mapsto \left[ F\left(\De_{[p\colon x\to y]\in P} p\sigma_{yx}\right) 
\colon \ce_X F(\B) \to \ce_Y F(\B) \right] , 
\end{gather*}
is an equivalence of categories. 
\end{reftheorem}

The quotient category $\Laby_n$ is formed by imposing (factoring out) certain relations; 
we refer to \cite{Fun}, Definition 23, for details. 
Essentially, the result of the procedure may be summarised thus. 
A maze is \emph{pure} if all passages carry the label $1$ (\emph{ibid.}, Definition 26). 
When typesetting such mazes, these labels will always be suppressed. 
By Theorem 9, \emph{ibid.}, a basis for $\Laby_n$ is provided by the pure mazes 
of cardinality at most $n$.  

One then has:  

\begin{reftheorem}[\cite{Fun}, Theorem 16] \label{S: Num}
The functor $\Phi$ above induces an equivalence of categories 
$$ 
\Phi\colon \Num_n \to \Lin(\Laby_n,\Mod).
$$
\end{reftheorem}

For our purposes, it shall be convenient to consider a skeletal version of $\Laby$, in 
which a single canonical set
$$ 
[k] = \{1,2,\dots,k\}
$$
is taken to represent each cardinality $k$. For example, we recall from \cite{Fun}, paragraph 5.1, that 
$\Laby_2$ has a skeleton consisting of the objects $\{\}$, $\{1\}$, and $\{1,2\}$, with a basis 
consisting of the mazes:
$$
\maze{ & }  , \
\maze{1 & \ar[l] 1} ,   \
\maze{1 & \upar[l] \downar[l] 1} ,   \
\maze{1 & \ar[l] \ar[dl] 1 \\ 2}  , \
\maze{1 & \ar[l]  1 \\ & \ar[ul] 2} , \
\maze{1 & \ar[l] 1 \\ 2 & \ar[l] 2} , \
\maze{1 & \ar[dl] 1 \\ 2 & \ar[ul] 2} .
$$

Throughout, 
we shall freely make use of the language of multi-sets, sets with possibly repeated elements. 
The \textbf{degree} 
or \textbf{multiplicity} of an element $x$, denoted $\deg x$, 
counts the number of repetitions. We agree to denote by 
$\#X$ the \textbf{support} of the multi-set $X$ --- that is: the underlying set --- 
and by $\abs{X}$ the \textbf{cardinality} --- by which is meant the number of elements counted 
with multiplicity. Thus 
$ \abs{X} = \sum_{x\in\#X} \deg_X x$.
For an elementary discussion of the theory of multi-sets, the reader is referred to Section 1 
of \cite{Fun}.

We shall be concerned exclusively with the classical algebraic functors $S^n$, $\Lambda^n$, and $\Gamma^n$. 
Their representations as labyrinth modules are shewn in Table \ref{T: Laby}, 
where we have let $\B^n=\gen{e_1,\dots,e_n}$.

\btab{cc}
\hline
\textbf{Functor} & \textbf{Labyrinth Module} \\
\hline
$S^n$ & $[k]\mapsto \gen{e^Z | \#Z=[k] \och \abs{Z}=n}$ \\
$\Lambda^n$ & $[k]\mapsto \gen{e^{\wedge Z} | \#Z=[k] \och \abs{Z}=n}$ \\
$\Gamma^n$ & $[k]\mapsto \gen{e^{[Z]} | \#Z=[k] \och \abs{Z}=n}$ \\
\hline
\etabc{The Classical Algebraic Functors as Labyrinth Modules.\label{T: Laby}}

\section{Comparison of Extensions}

One may, to a pair of 
 homogeneous functors over a binomial ring, each of degree $n$, associate at least four possible Ext-groups, 
depending upon the ambient category judged relevant: 
$\Func$, $\Pol_n$, $\Num_n$, or $\HPol_n$. Our first mission shall be 
to examine the interrelations betwixt these four types of extensions. 


\bth 			\label{S: ExtF} 
The category $\Pol_n$ is closed in $\Func$ under extensions, so that, in an exact sequence of functors as below, 
$G$ is polynomial if and only if $F$ and $H$ are:
$$ \xymatrix{0 \ar[r] & F \ar[r] & G \ar[r] & H \ar[r] & 0} $$
Consequently, $$ \Ext^1_{\Func}(H,F) = \Ext^1_{\Pol_n}(H,F)$$
when $F$ and $H$ are polynomial of degree $n$.

In particular, over $\Z$, there will be equalities
$$ \Ext^1_{\Func}(H,F) = \Ext^1_{\Pol_n}(H,F) = \Ext^1_{\Num_n}(H,F).$$
\eth

\bpr
Descension upon cross-effects produces an exact sequence of multi-functors: 
$$ \xymatrix{0 \ar[r] & \ce_{n+1} F \ar[r] & \ce_{n+1} G \ar[r] & \ce_{n+1} H \ar[r] & 0} $$
The middle functor vanishes if and only if the flanking ones do. 
\epr


%
%

\bex				\label{E: Num}
The category $\Num_n$ of numerical functors over an arbitrary binomial ring $\B$ 
is not closed under extensions. We give a counter-example to this effect. 

Let $\Z\binom x-$ be the free binomial ring in one variable (cf.~Theorem 5 of \cite{BR}), 
and let $\B$ be the quotient of 
$\Z\binom x-$ by the binomial ideal generated by $x^2$. 
Define the additive functor $G\colon{}_\B\XMod\to{}_\B\Mod$ by
\begin{gather*} 
\B^n\mapsto \B^n \\
\left[ P(x)\colon \B^m\to\B^n\right] \mapsto \left[ P(0)\colon \B^m\to\B^n\right].
\end{gather*}
It is not $\B$-linear, which one sees by considering the homomorphism $x\colon \B\to\B$, for the map 
$G(x)\colon \B\to\B $ is zero, while $ xG(1)= x \neq 0$.

However, there is 
an exact sequence: 
$$ 
\xymatrix{0\ar[r] & F \ar[r] & G \ar[r] & H \ar[r] & 0   } ,
$$
where the subfunctor $F(M)=xM$ and quotient functor $H(M)=M/xM$ are both linear. 
Multiplication by $x$ acts as zero on 
all three functors. 
\eex

By \emph{torsion}, here and in what follows, we shall always mean \emph{$\Z$-torsion}.

\blem
Let $\chi\colon R\to S$ be an homomorphism of unital rings. Suppose $\Im\chi$ is of finite index in $S$ and that 
$N$ is a torsion-free module over $S$. There is an isomorphism of $S$-modules:  
$$ 
\Hom_R(S,N) \cong N.
$$
\elem

\bpr
There is always an homomorphism of $S$-modules
$$ \phi\colon N \to \Hom_R(S,N), \quad y\mapsto y^\ast, $$
which is $S$-linear since 
$$ \phi(sy) = (sy)^\ast = y^\ast \circ s^\ast = \phi(y)\circ s^\ast = s\phi(y). $$
Moreover, we may always define  
$$ \psi\colon \Hom_R(S,N) \to N, \quad \alpha\mapsto \alpha(1),$$
and one has $\psi\circ\phi=1_N$,
although $(\phi\circ\psi)(\alpha)=\alpha(1)^\ast$ may not equal $\alpha$, 
nor will $\psi$ necessarily be $S$-linear, only $R$-linear:
$$
\psi(\chi(r)\alpha)=(\chi(r)\alpha)(1)=\alpha(\chi(r))=\chi(r)\alpha(1)=\chi(r)\psi(\alpha). 
$$

Under our assumption on finite index, however, 
we may, for any $s\in S$, find a positive integer $m$ and an $r\in R$ such that $ms=\chi(r)$. Then 
\begin{gather*}
m\alpha(s) = \alpha(ms)  = \alpha(\chi(r)) =  \chi(r)\alpha(1) =  ms\alpha(1) =  m\alpha(1)^\ast(s)   \\ 
m\psi(s\alpha) = \psi(ms\alpha) = \psi(\chi(r)\alpha) = \chi(r)\psi(\alpha) = ms\psi(\alpha),
\end{gather*}
and we may conclude using the dearth of torsion in $N$.
\epr 

We recall from \cite{Pelle}, Proposition 2.4, the existence of a small projective generator 
$$
Q_n = \GH ,
$$
inducing a category equivalence 
$$
\HPol_n \sim {}_{\GB}\Mod,
$$
where $\GB$ denotes the Schur algebra. In like wise, assuming $\B$ binomial, 
there will be a category equivalence 
$$
\Num_n \sim {}_{\BB}\Mod ,
$$
where $\BB$ represents an augmentation quotient of the free algebra on a matrix ring;
cf.~Theorem 5 of \cite{PF}. Theorem 22, \emph{ibid.}, asserts that 
the divided power map
$$
\gamma_n\colon \BB \to \GB, \quad [\sigma] \mapsto \sigma^{[n]} ,
$$
begets the forgetful functor $\HPol_n\to\Num_n$.

\bth \label{S: ExtFP} Assume $\B$ binomial.
Let $F$ be homogeneous and $G$ numerical of degree $n$. There is an exact sequence: 
$$ \xymatrixcolsep{.46pc} 
\xymatrix{0\ar[r] & \Ext^1_{\HPol_n}(F,\Nat(Q_n,G)) \ar[r] & \Ext^1_{\Num_n}(F,G) \ar[r] 
& \Nat_{\HPol_n}(F,\Ext^1_{\Num_n}(Q_n,G)) } $$
which, when $G$ is torsion-free, simplifies to: 
$$ \xymatrixcolsep{1.6pc}
\xymatrix{0\ar[r] & \Ext^1_{\HPol_n}(F,G) \ar[r] & \Ext^1_{\Num_n}(F,G) \ar[r] 
& \Nat_{\HPol_n}(F,\Ext_{\Num_n}^1(Q_n,G)) } $$
\eth

\bpr
The divided power map induces a spectral sequence (\emph{Base Change for $\Ext$})
$$ 
\Ext_{\HPol_n}^p(F,\Ext_{\Num_n}^q(Q_n,G)) \Rightarrow^p \Ext^{p+q}_{\Num_n}(F,G).
$$
The sequence enunciated in the theorem is simply the commencement of the associated 
five-term exact sequence. 

For the special case when $G$ is torsion-free, one applies the preceding lemma to shew 
 $\Nat(Q_n,G)\cong G$. The assumption that $\Im \gamma_n$ have finite index is a consequence of 
 Theorem 17 of \cite{PF}.
\epr

\section{Extension Theory for Numerical Functors}

In this section, we shall state and prove our main result, a formula for computing $\Ext^1$ 
of combinations of classical algebraic functors. We assume, from here and on, a binomial base ring $\B$.

We shall in the sequel make use of the following fact. 
Numerical functors over a $\Q$-algebra are semi-simple; the reason being that, in this situation, 
numerical and strict polynomial functors co-incide, and the latter functor category is 
well-known to be semi-simple. 
%

\bex 
\emph{Polynomial} functors over $\Q$-algebras, other than $\Q$ itself, might not be semi-simple. 
Modifying 
the construction of Example \ref{E: Num}, let $\B$ be the quotient of 
$\Q\binom x-$ by the binomial ideal generated by $x^2$; then, in fact, $\B\cong \Q[x]/(x^2)$.
As before, there will be an exact sequence of additive functors:
$$ 
\xymatrix{0\ar[r] & F \ar[r] & G \ar[r] & H \ar[r] & 0   } 
$$
It does not split, because the two flanking functors are linear, whereas the middle one is not. 
\eex

\blem \label{L: Torsion} 
Let $A$ be a torsion-free, unital ring and let $M$ and $N$ be finitely generated $A$-modules.
If $$\Ext^1_{\Q\otimes_\Z A}(\Q\otimes_\Z M,\Q\otimes_\Z N)=0,$$ then $\Ext^1_A(M,N)$ is a torsion group, hence finite.
\elem

\bpr
Consider an exact sequence: 
$$ \xymatrix{ 0\ar[r] & K\ar[r]^\beta & P \ar[r] & M \ar[r] & 0, } $$
with $P$ projective. 
Let $\theta\colon K\to N$. Suppose, towards a contradiction, that 
$n\theta$, for $n\in\Z^+$, never factors through $\beta$.

Since $P$ may be chosen finitely generated, any map $\Q\otimes_\Z P\to \Q\otimes_\Z N$ is of the form 
$q\eta$, with $\eta\colon P\to N$ and $q\in\Q$. Any factorisation of $\theta\colon\Q\otimes_\Z K\to \Q\otimes_\Z N$ 
through $\beta\colon \Q\otimes_\Z K\to \Q\otimes_\Z P$ would therefore yield a factorisation of 
$n\theta\colon K\to  N$ through $\beta\colon K\to P$, for some $n\in\Z^+$.
Consequently,  $\theta\colon\Q\otimes_\Z K\to \Q\otimes_\Z N$ does not factor through $\beta$, which means it is non-zero 
in 
$$
\Ext^1_{\Q\otimes_\Z A}(\Q\otimes_\Z M,\Q\otimes_\Z N).
$$ 
(Note that $\Q\otimes_\Z P$ is projective over $\Q\otimes_\Z A$.)

From this absurdity, one infers that some $n\theta$ must factor through $\beta$, so that $n\theta=0\in\Ext^1_A(M,N)$.
\epr 

A numerical functor is of \textbf{finite type} if it is finitely generated 
as a module over $\Laby_n$. 
In particular, it is of finite type when finitely generated over $\B$.

\bth 			\label{S: Torsion}
If $F$ and $G$ are numerical functors of finite type, each of degree $n$, then 
$\Ext^1_{\Num_n}(F,G)$ is finite.
\eth

\bpr
The category $\Laby_n$ is free, hence torsion-free (cf.~\cite{Fun}, Theorem 9). 
We made the remark above that all exact sequences of numerical functors split over a $\Q$-algebra.   
Lemma \ref{L: Torsion} will then apply.
\epr

\blem \label{L: Torsion Part} 
Let $M$ and $N$ be modules over a unital ring $A$.
Suppose $N$ is torsion-free and that $\Ext^1_A(M,N)$ is a torsion group. Consider a free resolution:
$$
\xymatrix{ R\ar[r]^\beta & Q\ar[r]^\alpha & P \ar[r] & M \ar[r] & 0 }
$$
The extension group $\Ext^1_A(M,N)$ will be isomorphic with the torsion subgroup of 
$$ \Coker\alpha^\ast = \Hom_A(Q,N)/\Im\alpha^\ast. $$
\elem

\bpr
Apply $\Hom_A(-,N)$ to the resolution: 
$$
\xymatrix{ \Hom_A(R,N) & \ar[l]_-{\beta^\ast} \Hom_A(Q,N) & \ar[l]_-{\alpha^\ast} \Hom_A(P,N) }
$$
Since the requested extension group is
$$
\Ext^1_A(M,N) = \Ker \beta^\ast / \Im \alpha^\ast,
$$ 
it will be sufficient to shew that 
$$\Ker\beta^\ast=\Set{\theta\in\Hom_A(Q,N) | \exists n\neq 0: n\theta\in\Im\alpha^\ast}.$$ 

If $n\theta=\eta\alpha\in\Im\alpha^\ast$, then 
$$
n\beta^\ast(\theta) = n\eta\alpha\beta =0\in\Hom_A(R,N),
$$
which is torsion-free, so $\beta^\ast(\theta)=0$. Consequently, $\theta\in\Ker\beta^\ast$.

Conversely, if $\theta\in\Ker\beta^\ast$ has $n\theta\notin\Im\alpha^\ast$ for all $n\neq 0$, then 
the homology class of $\theta$ represents an element of infinite order in $\Ext^1_A(M,N)$.
\epr

To prepare for our main theorem, we fix some notation. Let $e_1,\dots,e_n$ be the canonical basis of $\B^n$ 
and let $F$ be one of the 
classical algebraical functors $S^m$, $\Lambda^m$, $\Gamma^m$. 
To a multi-set $X$, with support included in $[m]$, 
we associate a monomial $e^{X}_F$ in a natural fashion: 
$$
e_F^X = \bca
\prod_{x\in X} e_x & \text{if $F=S^m$;} \\
\bigwedge_{x\in X} e_x & \text{if $F=\Lambda^m$;} \\ 
\prod_{x\in\#X} e_x^{[\deg x]} & \text{if $F=\Gamma^m$.}
\eca
$$
Using this notation, the functor $F$ will correspond to the labyrinth module 
$$ 
[k]\mapsto \gen{e^Z_F | \#Z=[k] \och \abs{Z}=m}.
$$
(Cf.~Table \ref{T: Laby}.)

Next, let us choose a system of generators for $F$. 
Suppose $B$ is a family of multi-sets $X$, fulfilling $\#X\subseteq[m]$ and $\abs{X}=m$, such that
the monomials
$$
\Set{e_F^X | X\in B}
$$
generate $F(\B^m)$ as a module over $\Laby_m$.

\bth 		\label{S: D}
Let $F\in\{S^m, \Lambda^m, \Gamma^m\}$, $G\in\{S^n, \Lambda^n, \Gamma^n\}$. 
 Construct a matrix $D$ according to the subsequent specifications. 
\begin{enumerate}
\item Index the rows of $D$  by pairs of multi-sets $(X,Y)$, where $X\in B$, $\abs{Y}=n$, and $\#Y=\#X$.
\item Index the columns of $D$ by pairs $((P_X)_{X\in B},Z)$, where $Z$ is a multi-set fulfilling 
$\#Z\subseteq [n]$ and $\abs{Z}=n$; and $(P_X)_{X\in B}$ belongs to a system of generators for 
the annihilator
$$
\Ann \sum_{X\in B} e_F^X \stackrel{\mathrm{def}}{=} \Set{(P_X)_{X\in B} | \sum_X P_X e_X^F = 0} \subseteq \Laby_{\max(m,n)}^B.
$$
\item Let $D_{(X,Y), (P,Z)}$ be the co-efficient of $e^Z_G$ in $P_X e^Y_G$.
\end{enumerate}
Then 
$$ 
\Ext^1_{\Num_{\max(m,n)}}(F,G)
$$ 
is isomorphic with the torsion subgroup of $\Coker D$.

In particular, when $F$ and $G$ are integral functors, 
the specified algorithm will compute $\Ext^1_\Func(F,G)$.
\eth

\bpr
Both $F$ and $G$ are left modules over the category $L=\Laby_{\max(m,n)}$, viewed as a unital ring. 
Given a multi-set $X$, let 
$L(\#X)$ denote the projective left $L$-module consisting of mazes with source $\#X$.
Note that 
$$
\Hom_L(L,G) \cong G = \gen{e^Z_G | {\abs{Z}}=n}
$$
and 
$$
\Hom_L(L(\#X),G) \cong I_{\#X}G = \gen{e^Y_G | \#Y=\#X\och \abs{Y}=n}.
$$
($I_{\#X}$ betokens the identity maze on $\#X$.)

By assumption, the $L$-module $F$ is generated by the monomials $e^X_F$, for $X\in B$, and so a  
 projective resolution is:
$$		\xymatrixcolsep{4pc} 
\xymatrix{ 
\displaystyle \bigoplus_{P\in \Ann \sum_{X\in B} e_F^X } L  \ar[r]^-{\sum_{P} P^\ast} & 
\displaystyle \bigoplus_{X\in B} L(\# X) \ar[r]^-{\sum_{X\in B} (e^X_F)^\ast} &
F \ar[r] & 0
} 
$$
All differentials are right multiplication by some element $w$ 
of the target module. 
Following standard, we have denoted such a map by $w^\ast$.

Apply the functor $\Hom_L(-,G)$: 
$$ 			\xymatrixcolsep{4pc} 
\xymatrix{
\displaystyle \bigoplus_{P\in \Ann \sum_{X\in B} e_F^X } \gen{e^Z_G |  \abs{Z}=n} 
& \ar[l]_-{\left(P_\ast\right)_{P}} 
\displaystyle \bigoplus_{X\in B} \gen{e^Y_G | \#Y=\#X \wedge \abs{Y}=n}
} 
$$
The group $\Ext^1_{\Num_{\max(m,n)}}(F,G)$ sought for is now, by Lemma \ref{L: Torsion Part} 
and Theorem~\ref{S: Torsion}, 
equal to the torsion subgroup of the co-kernel of this homomorphism, described by precisely the matrix $D$. 
\epr

In practice, this reduces the calculation of the Extension Group to the determination of 
the Smith Normal Form for a finite matrix $D$. 
Since one may restrict attention to pure mazes $P_X$, its entries shall be integers, so the normal 
form will always exist, even though
the base ring $\B$ may not be a principal ideal domain.

\section{Sample Computations}

For all of the following examples, we restrict attention to $\B=\Z$, as being the case of most common interest.
We write $\Ext^1$ as an abbreviation for $\Ext^1_{\Func}$.

\bex 
It will perhaps assist the reader if, as our first illustration of Theorem \ref{S: D},
 we select a calculation of already familiar extension groups. We thus consider the case
 of quadratic functors $F,G\in\{S^2, \Lambda^2, \Gamma^2\}$. 
It is required to compute $\Ext^1(F,G)$, and the results will, unsurprisingly, be found to conform with 
those of Table \ref{T: Splitting}. We remind the reader of the basis for $\Laby_2$ exhibited towards 
the end of \ref{A: Preliminaries}.

\emph{Case 1: $F=S^2$.} Since $e_1e_2$ generates $S^2$, we may choose $B=\{\{1,2\}\}$.  
The annihilator $\Ann e_1e_2 $
contains, apart from (multiples of) the element 
$$
\maze{1 & \ar[l] 1 \\ 2 & \ar[l] 2} -
\maze{1 & \ar[dl] 1 \\ 2 & \ar[ul] 2} ,
$$
all those mazes having domain $\{\}$ or $\{1\}$. Those of the latter kind may safely be ignored; they contribute 
merely zeroes to the matrix $D$. Its single row will be indexed by $Y=X=\{1,2\}$ and the matrix will 
take on an appearance as follows. 

For $G=S^2$: 
$$						
\begin{array}{c|c}
		&		
\maze{1 & \ar[l] 1 \\ 2 & \ar[l] 2} -
\maze{1 & \ar[dl] 1 \\ 2 & \ar[ul] 2}	 
\\ 
\hline 
e_1e_2 		&		0 
\end{array}
$$
leading to $\Ext^1(S^2,S^2)=0$. 

For $G=\Lambda^2$: 
$$						
\begin{array}{c|c}
					&		
\maze{1 & \ar[l] 1 \\ 2 & \ar[l] 2} -
\maze{1 & \ar[dl] 1 \\ 2 & \ar[ul] 2}	 
\\ 
\hline 
e_1 \wedge e_2 		&		2(e_1\wedge e_2) 
\end{array}
$$
leading to $\Ext^1(S^2,\Lambda^2)=\Z_2$. 

For $G=\Gamma^2$: 
$$				
\begin{array}{c|c}
			&		
\maze{1 & \ar[l] 1 \\ 2 & \ar[l] 2} -
\maze{1 & \ar[dl] 1 \\ 2 & \ar[ul] 2}	 
\\ 
\hline 
e_1e_2 		&		0 
\end{array}
$$
leading to $\Ext^1(S^2,\Gamma^2)=0$. 

\emph{Case 2: $F=\Lambda^2$.} Since $e_1\wedge e_2$ generates $\Lambda^2$, we may choose 
$B=\{\{1,2\}\}$. The annihilator $\Ann( e_1\wedge e_2) $
is, again excluding mazes with domain $\{\}$ and $\{1\}$, spanned by the two elements  
$$
\maze{1 & \ar[l] 1 \\ 2 & \ar[l] 2} +
\maze{1 & \ar[dl] 1 \\ 2 & \ar[ul] 2} , \qquad 
\maze{1 & \ar[l]  1 \\ & \ar[ul] 2} .
$$
The single row of $D$ will be indexed by $Y=X=\{1,2\}$ and the matrix will 
take on an appearance as follows. 

For $G=S^2$: 
$$					
\begin{array}{c|cc}		
			&		
\maze{1 & \ar[l] 1 \\ 2 & \ar[l] 2} +
\maze{1 & \ar[dl] 1 \\ 2 & \ar[ul] 2}	 
			&		
\maze{1 & \ar[l]  1 \\ & \ar[ul] 2}
\\ 
\hline 
e_1e_2 		&		2e_1e_2		&		e_1^2
\end{array}
$$
leading to $\Ext^1(\Lambda^2,S^2)=0$. 

For $G=\Lambda^2$: 
$$						
\begin{array}{c|cc}
					&		
\maze{1 & \ar[l] 1 \\ 2 & \ar[l] 2} +
\maze{1 & \ar[dl] 1 \\ 2 & \ar[ul] 2}	 
			&		
\maze{1 & \ar[l]  1 \\ & \ar[ul] 2}
\\ 
\hline 
e_1 \wedge e_2 		&		0			& 0  
\end{array}
$$
leading to $\Ext^1(\Lambda^2,\Lambda^2)=0$. 

For $G=\Gamma^2$: 
$$						
\begin{array}{c|cc}
					&		
\maze{1 & \ar[l] 1 \\ 2 & \ar[l] 2} +
\maze{1 & \ar[dl] 1 \\ 2 & \ar[ul] 2}	 
			&		
\maze{1 & \ar[l]  1 \\ & \ar[ul] 2}
\\ 
\hline 
e_1 e_2 		&		2e_1e_2			& 2e_1^{[2]}  
\end{array}
$$
leading to $\Ext^1(\Lambda^2,\Gamma^2)=\Z_2$. 

\emph{Case 3: $F=\Gamma^2$.} Since $e_1^{[2]}$ generates $\Gamma^2$, we may choose $B=\{\{1,1\}\}$.  
The annihilator $\Ann e_1^{[2]} $
contains, suppressing mazes having domain $\{\}$ or $\{1,2\}$, only (multiples of) the element 
$$
\maze{1 & \upar[l] \downar[l] 1} -
2\maze{1 & \ar[l] 1} , 
$$
The single row of $D$ will be indexed by $Y=X=\{1,1\}$ and the matrix will 
take on an appearance as follows. 

For $G=S^2$: 
$$						
\begin{array}{c|c}
			&		
\maze{1 & \upar[l] \downar[l] 1} -
2\maze{1 & \ar[l] 1} 	 
\\ 
\hline 
e_1^2 		&		0
\end{array}
$$
leading to $\Ext^1(\Gamma^2,S^2)=0$. 

For $G=\Lambda^2$: 
$$						
\begin{array}{c|c}
			&		
\maze{1 & \upar[l] \downar[l] 1} -
2\maze{1 & \ar[l] 1} 	 
\\ 
\hline 
0	 		&		0
\end{array}
$$
leading to $\Ext^1(\Gamma^2,\Lambda^2)=0$.  

For $G=\Gamma^2$: 
$$						
\begin{array}{c|c}
			&		
\maze{1 & \upar[l] \downar[l] 1} -
2\maze{1 & \ar[l] 1} 	 
\\ 
\hline 
e_1^{[2]}	&		0
\end{array}
$$
leading to $\Ext^1(\Gamma^2,\Gamma^2)=0$.  
\eex

We next seek to extend our calculations to higher degrees. These efforts are generally hampered by the  
vastly complicated structure of the category 
$\Laby_n$ for $n\geq 3$. 
The size and complexity of the matrix $D$ in Theorem \ref{S: D}, will expand accordingly. 
Fortunately, it shall be found unnecessary to determine the whole matrix in order to compute homology. 
The next lemma can be used with advantage in certain cases. 

\blem[The Pivot Lemma] Let $D=(d_{jk})\colon \Z^m\to\Z^n$ be an integral matrix satisfying:
\bnum
\item $d_{jk}=0$ if $j>k$. 
\item $d_{kk}\mid d_{jk}$ for all $j,k$. 
\enum
Then $$ \Coker D \cong \Z^{\max(n-m,0)} \oplus \bigoplus_{k} \Z_{d_{kk}}.$$
\elem

\bpr
Through elementary column operations, $D$ transforms into a matrix on diagonal form 
(though not necessarily square): 
$$
\begin{pmatrix}
d_{11} & 0 & \cdots \\
0 & d_{22} & \cdots \\ 
\vdots & \vdots & \ddots
\end{pmatrix}
$$ 
The procedure corresponds to a basis change in the source module $\Z^m$ and 
does not alter the co-kernel.
\epr

A matrix $D$, transferred to the form required by the theorem, is spoken of as being in 
\textbf{pivotal form}, and the diagonal elements $d_{kk}$ 
are referred to as \textbf{pivot elements}.
In the computation of the torsion subgroup of such a $D$, we need concern ourselves merely with 
the non-zero diagonal elements.

\bex
The group $ \Ext^1(I,S^n)$ was determined by Franjou in \cite{Puissances} 
(in fact, $\Ext^k(I,S^n)$ was successfully determined for any $k$).
We shew how to arrive at the result by using Theorem \ref{S: D}, combined with the Pivot Lemma. 

The apposite matrix $D$ has an especially simple form. 
Because $I$ is generated by the monomial $e_1$, the matrix has a single row, corresponding to 
the pair $(X,Y)=(\{1\},\{1^n\})$.  
The columns are indexed by mazes $P\colon \{1\}\to Z$ 
such that $Pe_1=0$, which holds whenever $\abs{P}\geq 2$. The entries of $D$ are the 
co-efficients from the expressions 
$Pe_1^n$, and the sought extension group is the torsion subgroup of $\Z_d$, where 
$d$ is the greatest common divisor of the elements of $D$.

It will be found sufficient, in this particular case, to restrict attention to \emph{simple} $P$; that is, those 
$P$ not containing any pair of parallel passages. For instance, in the case $n=4$, the mazes
$$
\xymatrixcolsep{1pc} \xymatrixrowsep{.7pc}
Q = 
\maze{1 & \ar[l] \ar[dl] \ar[ddl] 1 \\ 2 &  \\ 3 } \qquad 
R= \maze{1 & \upar[l] \downar[l] \ar[dl] 1 \\ 2 & }
$$
both kill $e_1$, while producing inside the matrix $D$: 
\begin{align*}
Qe_1^4 &= \binom{4}{2,1,1} e_1^2e_2e_3 + \binom{4}{1,2,1} e_1e_2^2e_3 + \binom{4}{1,1,2} e_1e_2e_3^2 \\ 
Re_1^4 &= \left( \binom{4}{2,1,1} + \binom{4}{1,2,1} \right) e_1^3e_2 + \binom{4}{1,1,2} e_1^2e_3^2 .
\end{align*}
Since we are merely seeking the greatest common divisor of the entries of $D$, the maze $R$ 
may be left out from consideration. 

With this simplification, the only mazes $P$ which need be considered are those of the form 
$$
\xymatrixcolsep{1pc} \xymatrixrowsep{.7pc}
P_k = 
\maze{1 & \ar[l] \ar[dl] \ar[dddl] 1 \\ 2  \\ \vdots \\ k}, \quad k\geq 2.
$$
One immediately finds 
$$ 
P_ke_1^n = \sum_{\substack{\sum m_i=n \\ m_i\geq 1}} \binom{n}{m_1,\dots,m_k} e_1^{m_1}\cdots e_k^{m_k} , 
$$
and so the number $d$ required is simply the greatest common divisor of all 
possible multinomial co-efficients 
$\binom{n}{m_1,\dots,m_k}$, for $k\geq 2$.

Bearing on the problem, we recall the following fact, well known to Number\hyp Theorists. 
Let $p$ denote a prime and let $n$ be a positive integer;
we enquire for the multiplicity of the divisor $p$ in $n!$.
The answer is contained in Legendre's Factorial Formula 
(\cite{Legendre}, paragraph \textsc{xvi} of the Introduction):
$$ 
\mu_p(n!)=\sum_{i=1}^\infty \left\lfloor \frac{n}{p^i} \right\rfloor,
$$
where $\mu_p(m)$ betokens the number of factors $p$ contained in the number $m$.
This formula we deploy for the computation of 
\begin{equation} \label{E: GCD}
d = \GCD_{k\geq 2} \binom{n}{m_1,\dots,m_k}
\end{equation}
as follows. 

\begin{itemize}
\item First observe that, since $\binom{n}{1,n-1}=n$, the number $d$ will divide $n$. 

\item If $n=1$, the greatest common divisor in \eref{E: GCD} is supposed to be taken over an empty 
set of multinomial co-efficients, which means $d=0$ in this case. 

\item If $n=p^j$ is a prime power, then 
$\binom{n}{m_1,\dots,m_k}$ is always divisible by $p$. This is well known 
for binomial co-efficients (it may easily be justified combinatorially), 
and the multinomial case follows since 
$$
\binom{p^j}{m_1,\dots,m_k}=\binom{p^j}{m_1}\binom{p^j-m_1}{m_2,\dots,m_k}.
$$  
On the other hand,
$$
\mu_p\binom{p^j}{p^{j-1},\dots,p^{j-1}} = \mu_p\left(\frac{p^j!}{(p^{j-1}!)^p}\right)
= \sum_{i=1}^\infty \left\lfloor \frac{p^j}{p^i} \right\rfloor 
- p\sum_{i=1}^\infty \left\lfloor \frac{p^{j-1}}{p^i} \right\rfloor = 1.
$$
Hence, the greatest common divisor $d=p$ in this case.

\item If $n=p^jn'$, for $n'>1$ and $\GCD(p,n')=1$, then 
\begin{align*}
\mu_p\binom{p^j n'}{p^j,p^j(n'-1)} &= \mu_p\left(\frac{(p^j n')!}{p^j!(p^j(n'-1))!}\right) \\
&= \sum_{i=1}^\infty \left( \left\lfloor \frac{p^j n'}{p^i} \right\rfloor 
- \left\lfloor \frac{p^j}{p^i} \right\rfloor 
- \left\lfloor \frac{p^j(n'-1)}{p^i} \right\rfloor\right) =0,
\end{align*}
since one readily verifies that the sum vanishes term-wise.
Hence the greatest common divisor $d=1$ in case $n$ has more than one prime factor. 
\end{itemize} 

We thus recoup the result from \cite{Puissances}: 
$$
\Ext^1(I,S^n) = 
\bca
\Z_p & \text{if $n=p^k$ is a prime power;}\\ 
0 & \text{if not.} \\

\eca
$$
\eex

\bex
The group $ \Ext^1(I,\Lambda^n)$ was likewise determined by Franjou in \cite{Puissances} 
(along with the higher $\Ext$-groups). We derive the result using Theorem \ref{S: D}. 

The matrix $D$ will have the same format as in the preceding example, but now its entries 
will be found by evaluating the co-efficients from 
$Pe_1^{\wedge n}$, where $P$ is any maze $\{1\}\to Z$ with $\abs{P}\geq 2$. 

But $e_1^{\wedge n}=0$ if $n\geq 2$, so $D$ will consist only of zeroes in 
this case, and the extension group will be trivial. 
We have thus re-established the result from \cite{Puissances}:  
$$
\Ext^1(I,\Lambda^n) = 0.
$$
(The case $n=1$ was covered by the preceding example.)
\eex

\bex \label{E: 3}
Our next goal is to calculate the group $\Ext^1(\Gamma^3,S^3)$. 
Since 
$$
\maze{1 & \ar[l] \ar[dl] 1 \\ 2 & \ar[dl] 2 \\ 3} 
e_1^{[2]}e_2 = e_1e_2e_3,
$$
the two monomials $e_1^{[3]}$ and $e_1^{[2]}e_2$ generate $\Gamma^3$, and we may choose 
$$ 
B=\{ \{1,1,1\}, \{1,1,2\}  \}. 
$$
The annihilator 
$$
\Ann ( e_1^{[3]} + e_1^{[2]}e_2 )  
$$
contains, among other things, the two mazes 
$$  
\maze{1 & \ar[l]  1 \\ & \ar[ul] 2}
- 3 \maze{1 & \ar[l] 1} 
, \qquad 
\maze{1 & \ar[l] 1 \\ 2 & \upar[l]\downar[l] 2},
$$
providing a matrix $D$ of the form:
$$		
\begin{array}{cc|ccc}
e^{[X]}		&e^{Y}&		
\maze{1 & \ar[l]  1 \\ & \ar[ul] 2}
- 3 \maze{1 & \ar[l] 1} 	 
			&		
\maze{1 & \ar[l] 1 \\ 2 & \upar[l]\downar[l] 2}\\ 
\hline 
e_1^{[3]} & e_1^3  
& -3e_1^3 & 0 & \cdots \\
e_1^{[2]}e_2^{[1]} & e_1^2e_2 
& e_1^3 & 0 & \cdots\\
& e_1e_2^2 
& e_1^3 & 2e_1e_2^2 & \cdots 
\end{array}
$$
Performing a change of basis leads to: 
$$		
\begin{array}{c|ccc}
&		
\maze{1 & \ar[l]  1 \\ & \ar[ul] 2}
- 3 \maze{1 & \ar[l] 1} 	 
			&		
\maze{1 & \ar[l] 1 \\ 2 & \upar[l]\downar[l] 2}\\ 
\hline  
e_1^{3} + 3e_1^{2}e_2 & 0 & 0 & \cdots \\
e_1^{2}e_2 & e_1^3 & 0 & \cdots\\
e_1 e_2^{2} - e_1^{2}e_2 & 0 & 2e_1e_2^2 & \cdots 
\end{array}
$$

We now proceed by the following observations.
\bnum
\item \emph{All entries of the first row are $0$.} This is clear, since an element 
of $\Ann ( e_1^{[3]} + e_1^{[2]}e_2 )$ obviously kills $e_1^{3} + 3e_1^{2}e_2$. 
\item \emph{All entries of the second row are divisible by $1$.} This is trivial.
\item \emph{All entries of the third row are divisible by $2$.}  
Consider an element $P\in \Ann(e_1^{[3]} + e_1^{[2]}e_2)$. We shall shew that the evaluation of 
$P(e_1 e_2^{2} - e_1^{2}e_2)$ contains only even co-efficients. 

Let us write 
\begin{align} \label{E: PQ}
P &= b \maze{1 & \ar[l] \ar[dl] 1 \\ 2}
+ a_{12} \maze{1 & \ar[l] 1 \\ 2 & \ar[l] 2} 
+ a_{21} \maze{2 & \ar[l] 1 \\ 1 & \ar[l] 2}
+ a_{13} \maze{1 & \ar[l] 1 \\ 3 & \ar[l] 2} \\
&\quad + a_{31} \maze{3 & \ar[l] 1 \\ 1 & \ar[l] 2}
+ a_{23} \maze{2 & \ar[l] 1 \\ 3 & \ar[l] 2}
+ a_{32} \maze{3 & \ar[l] 1 \\ 2 & \ar[l] 2} + Q, \nonumber
\end{align}
where $Q$ denotes a linear combination of mazes of types distinct from the ones listed. 
An inspection of Table \ref{T: Laby_3} gives at hand that these are the only ones that can possibly
provide \emph{odd} co-efficients before monomials of type $e_i^{[2]}e_j$, and so, recording only said monomials: 
\begin{align*}
0 = P(e_1^{[3]} + e_1^{[2]}e_2) \equiv \cdots & 
+ (b+a_{21})e_1 e_2^{[2]}    +(b+a_{12}) e_1^{[2]}e_2 \\
& + (b+a_{32})e_2 e_3^{[2]}    + (b+a_{23}) e_2^{[2]}e_3 \\
& + (b+a_{13})e_3 e_1^{[2]}    + (b+a_{31}) e_3^{[2]}e_1 \quad \mod 2.
\end{align*}
We conclude that all $a_{ji}\equiv -b\equiv a_{ij}$. 

Another glance at Table \ref{T: Laby_3} will convince us that the mazes listed in the formula 
\eref{E: PQ} are the only ones that produce odd co-efficients in the evaluation of 
$P(e_1 e_2^{2} - e_1^{2}e_2)$. It follows that 
\begin{align*}
P(e_1 e_2^{2} - e_1^{2}e_2) &\equiv 
(a_{12} - a_{21})e_1 e_2^{2}    +(a_{21}-a_{12}) e_1^{2}e_2 \\
& \quad + (a_{23} - a_{32})e_2 e_3^{2}    +(a_{32}-a_{23}) e_2^{2}e_3 \\
& \quad + (a_{31} - a_{13})e_3 e_1^{2}    +(a_{13}-a_{31}) e_3^{2}e_1 \equiv 0 \quad \mod 2,
\end{align*}
as desired. 
\enum

\btab{|cc|ccc|} 
\hline
&& $e_1^{[3]} + e_1^{[2]}e_2$ & $e_1e_2^2$ & $e_1^2e_2$   \\
\hline
\emph{\parbox{7.8em}{(Domain $\{1\}$)}} 
& $\maze{1 & \ar[l] 1}$  & $e_1^{[3]}$ & $0$ & $0$ \\
& $\maze{1 & \upar[l] \downar[l] 1}$  & $6e_1^{[3]}$ & $0$ & $0$ \\
& $\maze{1 & \upar[l] \ar[l] \downar[l] 1}$ & $6e_1^{[3]}$ & $0$ & $0$ \\ 
& $\maze{1 & \ar[l] \ar[dl] 1 \\ 2}$ & $e_1^{[2]}e_2+e_1e_2^{[2]}$ & $0$ & $0$ \\ 
& $\maze{1 & \upar[l] \downar[l] \ar[dl] 1 \\ 2}$ & $2e_1^{[2]}e_2$ & $0$ & $0$ \\ 
& $\maze{1 & \ar[l] \ar[dl] \ar[ddl] 1 \\ 2 \\ 3}$ & $e_1e_2e_3$ & $0$ & $0$ \\ 
\hline
\emph{\parbox{7.8em}{(Domain $\{1,2\}$, \\  Co-domain $\{1\}$)}}
& $\maze{1 & \ar[l]  1 \\ & \ar[ul] 2}$ & $3e_1^{[3]}$ & $e_1^3$ & $e_1^3$ \\
& $\maze{1 & \upar[l] \downar[l]  1 \\ & \ar[ul] 2}$ & $6e_1^{[3]}$ & $0$ & $2e_1^3$ \\
& $\maze{1 & \ar[l]  1 \\ & \upar[ul] \downar[ul] 2}$ & $0$ & $2e_1^3$  & $0$\\
\hline
\emph{\parbox{7.8em}{(Domain $\{1,2\}$, \\  Co-domain $\{1,2\}$)} }
& $\maze{1 & \ar[l] 1 \\ 2 & \ar[l] 2}$  & $e_1^{[2]}e_2$ & $e_1e_2^2$ & $e_1^2e_2$ \\ 
& $\maze{1 & \upar[l] \downar[l] 1 \\ 2 & \ar[l] 2}$ & $2e_1^{[2]}e_2$ & $0$ & $2e_1^2e_2$ \\ 
& $\maze{1 & \ar[l] 1 \\ 2 & \upar[l] \downar[l] 2}$ & $0$ & $2e_1e_2^2$ & $0$ \\
& $\maze{1 & \ar[l] \ar[dl] 1 \\ 2 & \ar[l] 2}$ & $2e_1^{[2]}e_2$ & $0$ & $2e_1^2e_2$ \\
& $\maze{1 & \ar[l] 1 \\ 2 & \ar[l] \ar[lu] 2}$ & $0$ & $2e_1e_2^2$ & $0$ \\
\hline
\emph{\parbox{7.8em}{(Domain $\{1,2\}$, \\ Co-domain $\{1,2,3\}$)} }
& $\maze{1 & \ar[l] \ar[dl] 1 \\ 2 & \ar[dl] 2 \\ 3 & }$ & $e_1e_2e_3$ & $0$ & $2e_1e_2e_3$ \\ 
& $\maze{1 & \ar[l] 1 \\ 2 & \ar[dl] \ar[l] 2 \\ 3 & }$ & $0$ & $2e_1e_2e_3$ & $0$ \\ 
\hline
\etabc{The Action of Pure Mazes of the Category $\Laby_3$. 
(The table comprehends all mazes with domain $\{1\}$ or $\{1,2\}$, up to left multiplication with 
a permutation.)
\label{T: Laby_3}}

We have thus established that the matrix $D$ found above is on pivotal form. 
The Pivot Lemma, in conjunction with Theorem \ref{S: D}, will then apply to shew that 
$$ \Ext^1(\Gamma^3,S^3) = \Z_2.$$
\eex

\bex \label{E: X} 
The elusive non-trivial element of $\Ext^1(\Gamma^3,S^3)$ deserves consideration, and we provide 
the details of its construction. Define the functor $X$ by 
$$
X(M) = \gen{p,\ x^{[3]},\ x^{[1]}y^{[1]}z^{[1]},\ x^{[2]}y^{[1]}+\frac{x^2y+xy^2}{2} | p\in S^3(M),\ x,y,z\in M} ,
$$
a subfunctor of $\left(\gen{\frac12}\otimes S^3\right)\oplus \Gamma^3$. 
(The symbol $\gen{Z}$ betokens the abelian group generated by $Z$.) 

There is an obvious exact sequence: 
$$ 
\xymatrix{0\ar[r] & S^3 \ar[r] & X \ar[r] & \Gamma^3 \ar[r] & 0   } 
$$
Here is an elementary argument as to why it does not split. 
Suppose there were a retraction $\rho\colon X\to S^3$, and consider its restriction to 
$\Z^2=\gen{e_1,e_2}$. The maze 
$$
P= \maze{1 & \ar[l] 1 \\ 2 & \upar[l]\downar[l] 2}
$$
has the following effect on symmetric and divided powers: 
$$
\left\{ \begin{array}{rcl}
e_1^{3} 			&\mapsto& 0 			\\ 
e_1^{2}e_2 		&\mapsto& 0			\\
e_1e_2^2 		&\mapsto& 2e_1e_2^2 	\\
e_2^{3} 			&\mapsto& 0			 			
\end{array} \right.
\qquad
\left\{\begin{array}{rcl}
e_1^{[3]} 				&\mapsto& 0 			\\
 e_1^{[2]}e_2^{[1]} 		&\mapsto& 0 			\\
e_1^{[1]}e_2^{[2]} 		&\mapsto& 2e_1^{[1]}e_2^{[2]} 	\\
e_2^{[3]} 				&\mapsto& 0.
\end{array} \right.
$$
An application of $P$ to a symmetric polynomial will thus
single out the term $e_1^{1}e_2^{2}$, doubling it. 

Since
$$ P\rho(2e_1^{[2]}e_2^{[1]}) = \rho(P\cdot 2e_1^{[2]}e_2^{[1]}) = \rho(0) = 0,$$
the co-efficient in $\rho(2e_1^{[2]}e_2^{[1]})\in S^3(\Z^2)$ of $e_1^{1}e_2^{2}$ must be $0$.
The same co-efficient in $2\rho\left(e_1^{[2]}e_2^{[1]}+\frac{e_1^2e_2+e_1e_2^2}{2}\right)$ 
must, at the very least, be even. We reach a contradiction upon inspection of
\begin{align*}
2\rho\left(e_1^{[2]}e_2^{[1]}+\frac{e_1^2e_2+e_1e_2^2}{2}\right) - \rho(2e_1^{[2]}e_2^{[1]})
&= \rho\left(e_1^2e_2+e_1e_2^2\right) = e_1^2e_2 + e_1e_2^2,
\end{align*}
in which expression the co-efficient of $e_1e_2^{2}$ is ostensibly odd. 
\eex

\bex				\label{E: 4}
A similar analysis may be carried out for the zenzizenzic case $\Ext^1(\Gamma^4,S^4)$.
Since 
\begin{gather*}
\maze{1 & \ar[l] \ar[dl] 1 \\ 2 } e_1^{[4]}
- \left( \maze{1 & \ar[l] 1 \\ 2 & \ar[l]  2  }
+ \maze{1 & \ar[dl] 1 \\ 2 & \ar[lu]  2  } \right) e_1^{[3]}e_2 = e_1^{[2]}e_2^{[2]}, \\
\maze{1 & \ar[l] 1 \\ 2 & \ar[l] \ar[dl] 2 \\ 3 } e_1^{[2]}e_2^{[2]} = e_1^{[2]}e_2e_3, \qquad
\maze{1 & \ar[l] \ar[dl] 1 \\ 2 & \ar[dl] \ar[ddl] 2 \\ 3 \\ 4 } e_1^{[2]}e_2^{[2]} = e_1e_2e_3e_4;
\end{gather*}
the monomials $e_1^{[4]}$ and $e_1^{[3]}e_2$ generate $\Gamma^4$, and we may choose 
$$ 
B=\{ \{1,1,1,1\}, \{1,1,1,2\}  \}. 
$$
The annihilator 
$$
\Ann ( e_1^{[4]} + e_1^{[3]}e_2 )  
$$
contains, among other things, the two mazes 
$$
\maze{1 & \ar[l]  1 \\ & \ar[ul] 2}
- 4 \maze{1 & \ar[l] 1} 
, \qquad 
\maze{1 & \ar[l] 1 \\ 2 & \ar[l] \ar[dl] 2 \\ 3 }, 
$$
providing a matrix $D$ of the form:
$$		
\begin{array}{cc|ccc}
e^{[X]}		& e^{Y}	
& \maze{1 & \ar[l]  1 \\ & \ar[ul] 2}
- 4 \maze{1 & \ar[l] 1}  	 		
& \maze{1 & \ar[l] 1 \\ 2 & \ar[l] \ar[dl] 2 \\ 3 }
\\ 
\hline 
e_1^{[4]} & e_1^4  & -4e_1^4 & 0 & \cdots \\ 
e_1^{[3]}e_2^{[1]} & e_1^3e_2 & e_1^4 & 0 & \cdots \\ 
& e_1^2e_2^2 & e_1^4 & 2e_1^2e_2e_3 & \cdots \\ 
& e_1e_2^3 & e_1^4  &  3e_1e_2^2e_3+3e_1e_2e_3^2 & \cdots \\ 
\end{array}
$$
A change of basis leads to: 
$$		
\begin{array}{c|ccc}
 &	\maze{1 & \ar[l]  1 \\ & \ar[ul] 2}
- 4 \maze{1 & \ar[l] 1}	 
& \maze{1 & \ar[l] 1 \\ 2 & \ar[l] \ar[dl] 2 \\ 3 } 
 \\ 
\hline 
 e_1^4+4e_1^3e_2  		& 0 & 0 &  \cdots \\
 e_1^3e_2 				& e_1^4 & 0 & \cdots \\ 
 e_1^2e_2^2-e_1^3e_2 		& 0 & 2e_1^2e_2e_3 & \cdots \\ 
 e_1e_2^3-e_1^3e_2 		& 0  &  3e_1e_2^2e_3+3e_1e_2e_3^2 & \cdots \\ 
\end{array}
$$

\btab{|c|ccc|} 
\hline
& $e_1^{[4]} + e_1^{[3]}e_2$ & $e_1^2e_2^2$ & $e_1^3e_2$   \\
\hline
$\maze{i & \ar[l] \ar[dl] 1 \\ j}$  
& $e_i^{[3]}e_j+e_i^{[2]}e_j^{[2]}+e_ie_j^{[3]}$ & $0$ & $0$ \\ 
 $\maze{i & \ar[l] \ar[dl] \ar[ddl] 1 \\ j \\ k}$  
& $e_i^{[2]}e_je_k+e_ie_j^{[2]}e_k+e_ie_je_k^{[2]}$ & $0$ & $0$ \\ 
\hline
$\maze{i & \ar[l]  1 \\ & \ar[ul] 2}$  
& $4e_i^{[4]}$ & $e_i^4$ & $e_i^4$ \\
\hline
$\maze{i & \ar[l] 1 \\ j & \ar[l] 2}$  
& $e_i^{[3]}e_j$ & $e_i^2e_j^2$ & $e_i^3e_j$ \\ 
 $\maze{i & \ar[l] \ar[dl] 1 \\ j & \ar[ul] 2}$  
& $2e_i^{[2]}e_j^{[2]}+3e_i^{[3]}e_j$ & $2e_i^3e_j$ & $3e_i^2e_j^2+3e_i^3e_j$ \\
\hline
$\maze{i & \ar[l] \ar[dl] 1 \\ j & \ar[dl] 2 \\ k & }$ 
& $e_i^{[2]}e_je_k+e_ie_j^{[2]}e_k$ & $2e_ie_je_k^2$ & $3e_i^2e_je_k+3e_ie_j^2e_k$ \\ 
\hline
\etabc{The Action of Pure Mazes of the Category $\Laby_4$.\label{T: Laby_4}}

We now proceed by the following observations.
\bnum
\item \emph{All entries of the first row are $0$.} This is clear, since an element 
of $\Ann ( e_1^{[4]} + e_1^{[3]}e_2 )$ must kill $e_1^4 + 4e_1^3e_2$. 

\item \emph{All entries of the second row are divisible by $1$.} This is trivial.

\item \emph{All entries of the third row are divisible by $2$.}  
Consider an element $P\in \Ann(e_1^{[4]} + e_1^{[3]}e_2)$. We shall shew that the evaluation of 
$P(e_1^2 e_2^2 - e_1^3 e_2)$ contains only even co-efficients. 

One may verify that the mazes catalogued in Table \ref{T: Laby_4} are the only ones 
that produce odd co-efficients before monomials of type $e_i^{[3]}e_j$, $e_i^{[2]}e_j^{[2]}$, or 
$e_i^{[2]}e_je_k$ in the evaluation of $P(e_1^{[4]} + e_1^{[3]}e_2)$. Let us, therefore, write 
\begin{align*} 
P &= 
\sum_{i,j} a_{i,j}\maze{i & \ar[l] 1 \\ j & \ar[l] 2}
+ \sum_{i,j} b_{i,j}\maze{i & \ar[l] \ar[dl] 1 \\ j & \ar[ul] 2}
+ \sum_{i,j} c_{ij}\maze{i & \ar[l] \ar[dl] 1 \\ j} \\
&\quad +\sum_{i}d_{i}\maze{i & \ar[l]  1 \\ & \ar[ul] 2}
+\sum_{i,j,k} p_{ij,k} \maze{i & \ar[l] \ar[dl] 1 \\ j & \ar[dl] 2 \\ k & } 
+ \sum_{i,j,k} q_{ijk} \maze{i & \ar[l] \ar[dl] \ar[ddl] 1 \\ j \\ k} 
+ Q ,
\end{align*}
where $Q$ denotes a linear combination of mazes of types distinct from the ones listed. 
Keeping track of monomials of type 
$e_i^{[3]}e_j$, $e_i^{[2]}e_j^{[2]}$, and $e_i^{[2]}e_je_k$ only, we find
\begin{align*} 
0&= P(e_1^{[4]} + e_1^{[3]}e_2) \\ 
& \equiv \dots + \sum_{i,j} a_{i,j} e_i^{[3]}e_j
+ \sum_{i,j} b_{i,j}\left(2e_i^{[2]}e_j^{[2]}+3e_i^{[3]}e_j\right) \\
&\qquad + \sum_{i,j} c_{ij}\left(e_i^{[3]}e_j+e_i^{[2]}e_j^{[2]}+e_ie_j^{[3]}\right) 
+\sum_{i,j,k} p_{ij,k} \left(e_i^{[2]}e_je_k+e_ie_j^{[2]}e_k\right) \\
&\qquad + \sum_{i,j,k} q_{ijk} \left(e_i^{[2]}e_je_k+e_ie_j^{[2]}e_k+e_ie_je_k^{[2]}\right) 
\quad \mod 2,
\end{align*} 
from which we conclude 
$$
\left\{ \begin{array}{rcl}
a_{i,j} + 3b_{i,j} + c_{ij} &\equiv& 0 \\
2b_{i,j} + c_{ij} &\equiv& 0 \\
p_{ij,k} + p_{ik,j} + q_{ijk} &\equiv& 0;
\end{array} \right.
$$
and hence 
$$
\bca
a_{i,j}\equiv b_{i,j} \\
c_{ij}\equiv 0 \\
p_{ij,k} \equiv p_{jk,i} \equiv p_{ki,j} \\
 q_{ijk} \equiv 0. 
\eca
$$ 

One next verifies that the mazes listed in Table \ref{T: Laby_4} 
are the only ones that will produce odd co-efficients in the evaluation of $P(e_1^2 e_2^2 - e_1^3 e_2)$.
It follows that 
\begin{align*}
P(e_1^2 e_2^2 - e_1^3 e_2) &\equiv 
\sum_{i,j} a_{i,j} \left(e_i^2e_j^2-e_i^3e_j\right)
+ \sum_{i,j} b_{i,j} \left(-e_i^3e_j-3e_i^2e_j^2\right) \\
&\quad +\sum_{i,j,k} p_{ij,k} \left(2e_ie_je_k^2-3e_i^2e_je_k-3e_ie_j^2e_k\right)
\equiv 0 \quad \mod 2,
\end{align*}
as desired. 

\item \emph{All entries of the fourth row are divisible by $3$.}  
Again, consider an element $P\in \Ann(e_1^{[4]} + e_1^{[3]}e_2)$. We shall shew that the evaluation of 
$P(e_1 e_2^3 - e_1^3 e_2)$ contains only co-efficients divisible by $3$. 

Write 
\begin{equation} 			\label{E: PR}
P = 
\sum_{i,j} a_{i,j}\maze{i & \ar[l] 1 \\ j & \ar[l] 2}
+ \sum_{i,j} c_{ij}\maze{i & \ar[l] \ar[dl] 1 \\ j} 
+\sum_{i}d_{i}\maze{i & \ar[l]  1 \\ & \ar[ul] 2} +R ,
\end{equation}
where $R$ denotes a linear combination of mazes of types distinct from the ones listed. 
One easily verifies that these are the only mazes that, in the 
evaluation of $P(e_1^{[4]} + e_1^{[3]}e_2)$, yield co-efficients not divisible by $3$ 
before monomials of type $e_i^{[3]}e_j$. We calculate 
$$
0 = P(e_1^{[4]} + e_1^{[3]}e_2)  
 \equiv \dots + \sum_{i,j} a_{i,j} e_i^{[3]}e_j 
+ \sum_{i,j} c_{ij}\left(e_i^{[3]}e_j+e_ie_j^{[3]}\right)  
\quad \mod 3,
$$
from which $a_{i,j}\equiv -c_{ij}\equiv a_{j,i}$. 

One next verifies that the mazes listed in \eref{E: PR} are the only ones that will
produce co-efficients not divisible by $3$ in $P(e_1 e_2^3 - e_1^3 e_2)$. Consequently, 
$$
P(e_1 e_2^3 - e_1^3 e_2) \equiv 
\sum_{i,j} a_{i,j} \left(e_ie_j^3-e_i^3e_j\right) \equiv 0 \quad \mod 3,
$$
as required.
\enum

We have thus found a pivot matrix $D$, and
the Pivot Lemma, in conjunction with Theorem \ref{S: D}, applies to shew that 
$$ \Ext^1(\Gamma^4,S^4) = \Z_2\oplus \Z_3.$$
\eex

\section{A Systematic Census} 

We now propose to conduct a systematic survey of Extension groups for any combination 
of classical algebraic functors covered by Theorem \ref{S: D}. 
Our argument will split up into nine cases, each case being indicated by an appropriate symbol. 

\subsubsection*{Case A: $\Ext^1(S^m,\Gamma^n)$.} 

Since $e_1\cdots e_m$ generates $S^m$, we may choose $B=\{[m]\}$. 
When $m>n$, there shall be no rows of $D$. 
When $m=n$, there is a single row, indexed by $X=Y=[m]$. This row stems from $P(e_1\cdots e_m)$, 
where $P\in \Ann e_1\cdots e_m$, and therefore only consists of zeroes. 

Suppose, then, that $m<n$. The rows of $D$ are indexed by pairs $(X,Y)$, fulfilling 
$\#Y=\#X=[m]$ and $\abs{Y}=n$. The maze 
$$ 
P_Y\colon [m]\to \Set{(p,j) | 1\leq p\leq m \och 1\leq j\leq \deg_Y p} \cong [n],
$$ 
defined by 
$$ 
P_Y=\Set{p\to (p,j) | 1\leq p\leq m \och 1\leq j\leq \deg_Y p},
$$
belongs to $\Ann e_1\cdots e_m$. Let $Z$ be another multi-set with $\#Z=\#X=[m]$ and $\abs{Z}=n$; then 
$$
P_Y e^{[Z]} = 
\bca 
\prod_{(p,j)} e_{(p,j)} & \text{if $Z=Y$;} \\
0 & \text{if $Z\neq Y$.} 
\eca
$$
By the Pivot Lemma, the co-kernel is then trivial, and we have shewn that 
$$ 
\Ext^1(S^m,\Gamma^n)=0 \quad \text{for all $m,n$.}
$$

\subsubsection*{Case B: $\Ext^1(S^m,\Lambda^n)$.} 

Since $e_1\cdots e_m$ generates $S^m$, we may choose $B=\{[m]\}$. 
When $m>n$, there is no monomial $e^{\wedge Y}$ with $\#Y=\#X=[m]$ and $\abs{Y}=n$. When $m<n$, all such monomials 
will contain repetitions and so be zero. Only in the case $m=n$ does $D$ contain a positive number of rows. 
Its single row is indexed by $X=Y=[m]$, and the entries are given by the co-efficients of 
$P(e_1\wedge\dots\wedge e_m)$, where $P\in \Ann e_1\cdots e_m$. These are easily seen to be even when 
$m=n\geq 2$, and so the Pivot Lemma gives 
$$
\Ext^1(S^m,\Lambda^n)= 
\bca
\Z_2 & \text{when $m= n\geq 2$;} \\
0 & \text{otherwise.}
\eca 
$$

\subsubsection*{Case C: $\Ext^1(S^m,S^n)$.} 

Since $e_1\cdots e_m$ generates $S^m$, we may choose $B=\{[m]\}$. When $m> n$, there shall be no rows of $D$. 
When $m=n$, there is a single row, indexed by $X=Y=[m]$. This row stems from $P(e_1\cdots e_m)$, 
where $P\in \Ann e_1\cdots e_m$, and therefore only consists of zeroes. 
We have thus shewn that 
$$ 
\Ext^1(S^m,S^n)=0 \quad \text{when $m\geq n$.}
$$

\subsubsection*{Case D: $\Ext^1(\Lambda^m,\Gamma^n)$.} 

Since $e_1\wedge\cdots\wedge e_m$ generates $\Lambda^m$, we may choose $B=\{[m]\}$. 
For $m>n$, there shall be no rows of $D$. 
For $m<n$, one may again consider the mazes $P_Y$, defined in the discussion pertaining to Case A. 
Each step of the process there carried out may be repeated without alteration. 
When $m=n$, finally, there is a single row of $D$, indexed by $X=Y=[m]$. The entries are given 
by the co-efficients of 
$P(e_1\cdots e_m)$, where $P\in \Ann (e_1\wedge\cdots\wedge e_m)$. These are easily seen to be even when 
$m=n\geq 2$, and so the Pivot Lemma gives 
$$
\Ext^1(\Lambda^m,\Gamma^n)= 
\bca
\Z_2 & \text{when $m=n\geq 2$;} \\
0 & \text{otherwise.}
\eca 
$$

\subsubsection*{Case E: $\Ext^1(\Lambda^m,\Lambda^n)$.} 

Since $e_1\wedge\cdots\wedge e_m$ generates $\Lambda^m$, we may choose $B=\{[m]\}$. 
For $m>n$, there is no monomial $e^{\wedge Y}$ with $\#Y=\#X=[m]$ and $\abs{Y}=n$. For $m<n$, all such monomials 
will contain repetitions and so be zero. When $m=n$, finally, there is a single row, indexed by  
 by $X=Y=[m]$. 
 This row stems from $P(e_1\wedge\cdots\wedge e_m)$, 
where $P\in \Ann (e_1\wedge\cdots\wedge e_m)$, and therefore only consists of zeroes. 
We have thus shewn
$$ 
\Ext^1(\Lambda^m,\Lambda^n)=0 \quad \text{for all $m,n$.}
$$

\subsubsection*{Case F: $\Ext^1(\Lambda^m,S^n)$.} 

Since $e_1\wedge\cdots\wedge e_m$ generates $\Lambda^m$, we may choose $B=\{[m]\}$. 
When $m>n$, there is no monomial $e^{Y}$ with $\#Y=\#X=[m]$ and $\abs{Y}=n$. 
When $m=n$, the matrix $D$ contains a single row, indexed by $X=Y=[m]$, and the entries are 
given by the co-efficients of 
$P(e_1\cdots e_m)$, where $P\in \Ann (e_1\wedge\cdots\wedge e_m)$. 
Since 
$$ 
P=\Set{p\to 1 | 1\leq p\leq m}\colon [m] \to [1]
$$
is such a maze, mapping 
$ P(e_1\cdots e_m) = e_1^m$,
the Pivot Lemma shews that the co-kernel is trivial. 
We have thus proved
$$
\Ext^1(\Lambda^m,S^n)= 0 \quad \text{when $m\geq n$.} 
$$

\subsubsection*{Case G: $\Ext^1(\Gamma^m,\Gamma^n)$.} 

We consider the case $m=n$ only.
Choose for $B$ one representative of each isomorphism class of multi-sets of cardinality $m$, among which we may assume $\{1^m\}$.
When constructing the matrix $D$, we effectuate a change of basis as follows. Add the rows corresponding to 
all $Y=X\in B$ to the row corresponding to 
$Y=X=\{1^m\}$. The latter row shall thus set out the co-efficients of $P\left(\sum_{Z\in B}e^{[Z]}\right)$, 
for $P\in\Ann\sum_{Z\in B}e^{[Z]}$, and hence consist exclusively of zeroes. 

The remaining rows are indexed by monomials $(X,Z)$ such that $\#Z=\#X=[m]$ and $\abs{Z}=n$.
Define the mazes $P_Y$ as in Case A. For $Y\in B$, we find 
$$
\Ann\sum_{Z\in B}e^{[Z]} \ni P_Y - P_{\{1^m\}} \colon 
e^{[Z]} \mapsto 
\bca
e_1\cdots e_m & \text{if $Z=Y$;} \\
0 & \text{if $Z\neq Y$;}
\eca 
$$ 
and for $Y\notin B$
$$
\Ann\sum_{Z\in B}e^{[Z]} \ni P_Y \colon 
e^{[Z]} \mapsto 
\bca
e_1\cdots e_m & \text{if $Z=Y$;} \\
0 & \text{if $Z\neq Y$.}
\eca 
$$
The Pivot Lemma is then brought to apply, shewing 
$$ 
\Ext^1(\Gamma^m,\Gamma^n) = 0 \quad \text{when $m=n$.}
$$

\subsubsection*{Case H: $\Ext^1(\Gamma^m,\Lambda^n)$.} 

Choose for $B$ one representative $X$ of each isomorphism class of multi-sets of cardinality $m$, among which we may assume $\{1^m\}$.
When $m<n$, all monomials $e^{\wedge Y}$ with $\#Y=\#X\in B$ and $\abs{Y}=n$ will contain repetitions and so be zero. 
When $m=n$, the matrix $D$ contains a single row, indexed by $X=Y=[m]$, and the entries are given by the co-efficients of 
$P(e_1\wedge\cdots\wedge e_m)$, where $P\in \Ann \sum_{Z\in B} e^{[Z]}$.
Since 
$$ 
\Ann\sum_{Z\in B}e^{[Z]} \ni P_{[m]} - P_{\{1^m\}}  \colon e_1\wedge\cdots\wedge e_m \mapsto e_1\wedge\cdots\wedge e_m,
$$
the Pivot Lemma shews that the co-kernel is trivial. 
We have thus proved
$$
\Ext^1(\Gamma^m,\Lambda^n)= 0 \quad \text{when $m\leq n$.} 
$$

\subsubsection*{Case I: $\Ext^1(\Gamma^m,S^n)$.} 

We have found no conclusions of a general nature pertaining to this case, 
save for the calculations in low degrees 
performed in Example \ref{E: 3}. This would appear the most troublesome case. 

\medskip
We summarise our findings: 

\begin{theorem} \label{S: Ext}
Over the base ring $\Z$, one has the following extension groups in the functor category $\Func$:
\balph
\item $\Ext^1(S^m,\Gamma^n)=0$ for all $m,n$.

\item $\Ext^1(S^m,\Lambda^n)= 
\bca
\Z_2 & \text{when $m=n\geq 2$;} \\
0 & \text{otherwise.}
\eca $

\item $ \Ext^1(S^m,S^n)=0$ when $m\geq n$.

\item $\Ext^1(\Lambda^m,\Gamma^n)= 
\bca
\Z_2 & \text{when $m=n\geq 2$;} \\
0 & \text{otherwise.}
\eca $ 

\item $\Ext^1(\Lambda^m,\Lambda^n)=0$ for all $m,n$.

\item $\Ext^1(\Lambda^m,S^n)= 0$ when $m\geq n$. 

\item $\Ext^1(\Gamma^m,\Gamma^n) = 0$ when $m=n$.

\item $\Ext^1(\Gamma^m,\Lambda^n)= 0$ when $m\leq n$. 

\ealph \hfill
\end{theorem}

\end{document}